\documentclass[12pt]{amsart}
\usepackage{amsfonts,amsmath,amscd,amssymb}
\numberwithin{equation}{section}
\newtheorem{thm}{Theorem}[section]

\newtheorem{lem}[thm]{Lemma}
\newtheorem{defn}[thm]{Definition}
\newtheorem{cor}[thm]{Corollary}
\newtheorem{rem}[thm]{Remark}
\def\Lie{\mathop{\operatorname{\mathrm Lie}}\nolimits}
\def\Ind{\mathop{\operatorname{\mathrm Ind}}\nolimits}
\def\ind{\mathop{\operatorname{\mathrm Ind}}\nolimits}
\def\End{\mathop{\operatorname{\mathrm End}}\nolimits}
\def\vol{\mathop{\operatorname{\mathrm vol}}\nolimits}
\def\tr{\mathop{\operatorname{\mathrm tr}}\nolimits}
\def\Ad{\mathop{\operatorname{\mathrm Ad}}\nolimits}
\begin{document}
\title{Quantization of Fields and Automorphic Representations}
\author{Do Ngoc Diep}
\maketitle
\begin{abstract}
In this paper we use the quantization of fields based on Geometric Langlands Correspondence \cite{diep1} to realize the automorphic representations of some concrete series of groups: for the affine Heisenberg (loop) groups it is reduced to the construction of the affine Kac-Moody representation by the Weyl relations in Fock spaces. For the solvable and nilpotent groups following the construction we show that it is the result of applying the constructions of irreducible unitary representation via the geometric quantization and the construction of positive energy representations  and finally, for the semi-simple or reductive Lie groups, using the Geometric Langlands Correspondence, we show that a repeated application of the construction give all the automorphic representations of reductive Lie groups: first we show that every representation of the fundamental group of Riemann surface into the dual Langlands groups ${}^LG$ of $G$ corresponds to a representation of the fundamental group of the surface into the reductive group $G$, what is corresponding to a quantum inducing bundle of the geometric quantization of finite dimensional reductive Lie groups and  then apply the construction of positive energy representation of loop groups.
\end{abstract}
{\it MSC2010 Subject Classification: 

Primary: 22E57; Secondary:11F70, 14D24, 53D50}


\section{Introduction}
In the previous work \cite{diep1} we introduced a new quantization procedure based on the Geometric Langlands Correspondence. More precisely, started from fields in the target space, we proposed to use reduction to the case of fields on one complex variable target space, which is a Riemann surface, by using the the reduction of Kaluza-Klein of the extra dimension and symmetry, 
as the dual Langlands group ${}^LG$. Use the sigma model we can reduce it to a problem of quantization of the trivial vector bundle with connections over the space dual to the Lie algebra ${}^L\mathfrak g = \Lie({}^LG)$. By using the ideas of electric-magnetic duality we pass to the same problem for the Langlands dual group $G$. We then have a representation of the fundamental groups of the Riemann surface base. Then, we have some affine Kac-Moody loop algebra of meromorphic functions with value in Lie algebra $\mathfrak g=\Lie(G)$. 

Our procedure of quantization consists of two steps:
\begin{itemize}
\item Quantize the Lie algebra $\mathfrak g$ by using the Orbit Method to construct irreducible representations in the space $\mathcal H = \ind_P^G L^2(G; \mathfrak p,F,\rho,\sigma)$.
\item Construct the Fock space $\exp \mathcal H$ and the corresponding representations of the affine Kac-Moody loop algebra in it.
\end{itemize}

The main goal of this paper is to realize these steps for concrete series of Lie groups. In Section 2 we compute the representations of the affine Heisenberg groups. It is precisely the affine Weyl representations in the Fock space.
In Section 3 and 4 we discribe the results for the affine nilpotent or solvable Kac-Moody loop algebras and the corresponding groups. 
In the last section 5, we study the case of affine reductive or semi-simple Lie loop algebras and groups.
The main result (Theorem \ref{thm51} ) is obtained in Section 5. We show that every authomorphic representation can be obtained from our construction by decomposing the parabolic subgroup into a tower of parabolic subgroups, which are maximal in each-another.

\section{Irreducible Representations of Affine Heisenberg Groups}

\subsection{Orbit Method for Heisenberg Lie groups}
It is well-known that the Heisenberg algebra $\mathfrak g= \mathfrak h_{2n+1}$ is generated by $2n+1$ generators $X_1,\dots,X_n,Y_1,\dots, Y_n, Z$ satisfying the Weyl commutation relations 
\begin{equation}[X_i,Y_j] = Z\delta_{ij}I,   \end{equation}
Because the center of $\mathfrak h_{2n+1} $ is generated by $Z$, under an irreducible representation, it should be acted as some scalar, and we have the 
Weyl representations

\begin{equation} 
\left\{\begin{array}{rcl}
Z & \mapsto & \lambda\\
Y_j & \mapsto & x_j \\
X_i & \mapsto & \lambda\frac{\partial}{\partial x_i}, \lambda\in \mathbb R
\end{array}\right.
\end{equation}

In the particular case $\lambda= \frac{i}{\hbar}, i=\sqrt{-1}$, we have the ordinary Weyl commutation relations.

\subsection{Construction of
highest weight positive energy representations of loop algebras $\hat{\mathfrak g}$}

For any loop $T\in \hat{\mathfrak g},$ we define the generators with values in the representation of $\mathfrak g$,
$$L_n = \frac{1}{2\pi i}\oint_C z^{n+1}T(z)dz, \bar{L}_m = \frac{1}{2\pi i} \oint_Cdz \bar{z}^{n+1} \bar{T}(z).$$
We obtain then the
Virasoro relations:
$$[L_n,L_m] = (n-m)L_{m+n} + \frac{c}{12}n(n^2-1)\delta_{n+m,0}$$
$$[L_n,\overline{L}_m] = 0$$
$$[\overline{L}_n,\overline{L}_m] = (n-m)\overline{L}_{n+m} + \frac{c}{12}n(n^2-1)\delta_{n+m,0}$$
This relation let us to use the construction of representation of positive energy of loop algebras in Fock spaces.

\subsection{Fock space representations} 
Introduce the {\it vacuum state} $|0\rangle$ which is annihilated by all $L_n, n\geq 0$.
The highest state $|h\rangle$ is defined as vector such that
$$L_0|h\rangle = h|h\rangle, $$ and for all $n > 0$
$$L_0(L_{-k_1}\dots L_{-k_n}|h\rangle = (h+l)L_{-k_1} \dots L_{-k_n}|h\rangle$$ where $l = \sum_{i=1}^n k_i (k_i > 0).$

\begin{defn}[Verma module]
$$V_{c,h} = \{L_{k_1}\dots L_{-k_n}|h\rangle ; 0 \leq k_1 \leq \dots \leq k_n, n\in \mathbb N \}$$ of central charge $c$ and highest weight $h$
$$\bar{V}_{c,h} = \{\bar L_{-k_1}\dots \bar L_{-k_n}|h\rangle; 0 \leq k_1 \leq \dots \leq k_n, n\in \mathbb N  \}$$ of central charge $c$ and highest weight $\bar h$
\end{defn}
\begin{rem} The {\it total state space} is the direct sum
$$\sum_{h,\bar h}V_{c,h} \otimes \bar{V}_{c,\bar h}$$

The {\it singular vector} $v$ such that $L_n|v\rangle = 0, \forall n >0$. By quotienting out of $V_{c,h}$, the null submodule, we obtain an irreducible component representation $M_{c,h}$. 
\end{rem}

Then the quantum state space is $$\mathcal H = \bigoplus_{c,h} n_{c,h}M_{c,h} \otimes \bar M_{c,h}.$$

\begin{cor}
The obtained representation of $\hat{\mathfrak g} = L(C,\Ind(\mathfrak p,F,\sigma)(\mathfrak g))$ in the corresponding Fock space $\mathcal H$ is automorphic.
\end{cor}

By the construction of Fock space, it is clear that
\begin{lem}
The obtained representation is a subrepresentation of the exponential representation.
\end{lem}

Indeed there is a natural map from the exponential representation to the space of Fock representation of loop algebras.
We can therefore state the obtained result in the following form

\subsection{Construction of loop algebra representations}
The construction is reduced to the Fock space representations of the corresponding Kac-Moody Heibenberg algebras.
The construction is well-known and we refer the readers to \cite{keegan}
Now, affine Lie algebra of loop on $C$ 
$$\hat{\mathfrak g}=L(C,\Ind(\mathfrak p.F,\sigma)(\mathfrak g))$$ 
with values in the image
$$\Ind(\mathfrak p.F,\sigma)(\mathfrak g) \subseteq \End(\mathcal H)$$ 

\subsection{Representations of affine Heisenberg algebras}
Now use use two steps together then we have irreducible unitary representations of the affine Heisenberg groups in the corresponding Fock space.

Every element of the affine Kac-Moody Heisenberg algebra $\hat{\mathfrak g}$ is a holomorphic function on the Riemann surface with coefficients in our Heisenberg Lie algebra
 \begin{equation} X = \sum_{-\infty}^{+\infty}c_n z^n, c_n \in \mathfrak g \end{equation}
After the first step of quantization, $c_n$ become (un-)bounded but closed operators in the the space $\mathcal H = L^2(\mathbb R^k)$ 

We consider the tensor product of the ordinary affine Kac-Moody loop algebra and algebra $\hat{\mathfrak g}$.  
We have therefore,

\begin{thm}
Every irreducible unitarizable representation of the affine Kac-Moody Heisenberg groups can be obtained in this way, as the compositon of geometric quantization for finite dimensional Heisenberg algebras and the construction of positive energy representations of loop algebras.
\end{thm}
{\sc Proof.} Because we fixed a affine connection, we can decompose the element $ S\in \hat{\mathfrak g}$ as
\begin{equation} X = \sum_i X_i \otimes f_i, f_i \in Loop(\mathbb C) \end{equation} The sum is finite and the elements $X_i$ are operators in the space $\mathcal H = L^2(\mathbb R^k).$
 \hfill$\Box$

\section{Irreducible Representations of Affine Nilpotent Lie Groups}

The representations of nilpotent Lie groups are obtained by choosing a real or complex positive polarization and reduce to the Weyl commutation relations.

\begin{thm}
Every irreducible unitarizable representation of the affine Kac-Moody groups of connected nilpotent groups can be obtained in this way, as the composition of geometric quantization for finite dimensional Heisenberg algebras and the construction of positive energy representations of loop algebras.
\end{thm}
{\sc Proof.} 
Let us denote by $\mathfrak p\subset \mathfrak g_\mathbb{C}$ a positive polarization at $F \in \mathfrak g^*$. Then there is a natural map from $\mathfrak g \twoheadrightarrow \mathfrak p$ The Weyl representation of $\mathfrak p / \ker F|_{\mathfrak p}$ gives rise to the representation from the orbit method. 
\hfill$\Box$

\section{Irreducible Representations of Affine Solvabble Lie Groups}
The representations of solvable Lie groups are also obtained by choosing a real or complex positive polarization and reduce to the Weyl commutation relations.

\begin{thm}
Every irreducible unitarizable representation of the affine Kac-Moody solvable groups can be obtained in this way, as the compositon of geometric quantization for finite dimensional Heisenberg algebras and the construction of positive energy representations of loop algebras.
\end{thm}
{\sc Proof.} 
For a positive polarization $\mathfrak p / \ker F|_{\mathfrak p}$ is the Heisenberg Lie algebra and therefore one use the orbit method to iterate to the representation of the Lie algebra $\mathfrak p$ and $\mathfrak g$, in the space 
$\mathcal H = L^2(G; \mathfrak p, F, \rho, \sigma)$ 
of partially holomorphic partially invariant sections of the inducing bundle,  see \cite{kirillov} for more detailed expositon. Now apply the construction of Fock space 
\begin{equation}
\exp \mathcal H = \varinjlim_{n\to \infty} \sum_{k=0}^n \frac{1}{k!}\mathcal H ^{\otimes k}\end{equation}
\hfill$\Box$

\section{Automorphic Representations of Affine Reductive Lie Groups}
The automorphic representations are obtained by holomorphic induction and can be realized as the discrete series representation from the orbit method.

\subsection{Spectral side of the Arthur's trace formula}
We refer the readers to the recent work of T. Finis, E. M. Lapid and W. M\"uller \cite{flm} for the spectral side of the Arthur's trace formula. For simplicity, we do not use the globalized language of ad\`ele and id\`ele.

Let us consider a reductive group $G=\underline{G}(\mathbb R)$ defined over a number field $F$, for simplicity, often referred to as the field of rational numbers $\mathbb Q$.
Fix a maximal torus $T_0 \subset G$ and a maximal compact subgroup $K\subset G$. Denote $\mathcal P$ the set of parabolic subgroups of $G$ containing $T_0$. 
For each parabolic subgroup $P \in \mathcal P$, consider the Levi decomposition 
$P = M.N$ where $N = N(P)$ is the nilpotent radical of $P$. and $M$ is the Levi component of $P$. Denote $T_M$ the split part of the center of $M$ and $A_M=T_M(\mathbb R)^0$ the corresponding split torus. Its Lie algebra $\mathfrak a_M = \Lie A_M$ is generated by the lattice of co-roots of $T_M \cap (G,G)$. The dual space of $\mathfrak a_M$ is denoted by $(\mathfrak a_M^{G})^*$, and $r= \dim\mathfrak a_M^{G}$. Denote $N_G(M)$ the normalizer of $M$ in $G$ and $W(M) = N_G(M)/M$ the corresponding Weyl group. For any element $s\in W(M)$, $M_s = \langle M,s\rangle$ is the smallest Levi subgroup containing $M$ and $s$. Denote $\mathcal P(M)$ the set of all parabolic subgroup of $G$ containing $M$.

Fixed a parabolic subgroup $P$ denote $\Delta^G_P$ the set of all simple root and $\Sigma^G_P$ the set of simple reduced coroots of $P$ and $\bar{P}$ the opposite parabolic.

The elements of $\Delta^G_P$ could be considered as some $r$-tuple with a prescribed ordering introduced by a Weyl chamber.

Two parabolic subgroups $P$ and $Q$ are adjacent along $\alpha^chec$ iff 
\begin{equation}\Sigma^G_P \cap \Sigma^G_Q = \{\check{\alpha} \} \end{equation} 
i.e. $P$ and $Q$ are maximal in their product group $PQ$. The elements $\mu_1, \dots, \mu_r \in (\mathfrak a^G_M)^*$ are in general position, if 
for linear independent $\check{\alpha}, \dots, \check{\alpha}_r\in \Sigma^G_M$ and the dual basis $w_1,\dots,w_r\in (\mathfrak a^G_M)^*$ the parabolics $P_i$ and $Q_i$ are adjacent iff the element
\begin{equation}\sum_{j\ne i} \langle \mu_j,\check{\alpha}_j \rangle w_j \mbox{ lies in the Weyl chamber of } P_iQ_i \end{equation}

For $P\in \mathcal P(M)$ consider the induced representation
\begin{equation} \label{ind}
\mathcal A^2_P = \Ind_P^G L^2_{disc}(A_M\Gamma_M \setminus M)
  \end{equation}

When $\lambda$ is changing in $i(\mathfrak a^G_M)^*$ there
is a family of sub-quotients $I_P(\lambda)$. The theory of Eisenstein series gives rise to intertwining operators from $I_P(\lambda)$ to the authomorphic representations of $G$. 

The main result of T. Finis, E. M. Lapid and W. M\"uller (Theorem 1 in \cite{flm}) says that {\it the spectral side of Arthur's trace formula for decomposition of $L^2(\Gamma\setminus G)$ is given by the sum-integal 
of the expression

\begin{equation}
\vol(\mathfrak a^G_M/L)\frac{(-1)^k}{k!} |\det(x-1|\mathfrak a^{M_s}_M)  |^{-1} \times \end{equation}
$$ \tr(M_{P|Q_k}(\lambda)M'_{P_k|Q_k}(\langle \lambda,\alpha'_k\rangle)M_{P_k|Q_{k-1}}(\lambda) \dots M_{P_2|Q_1}(\lambda)$$ 
$$ M'_{Q_1|P_1}(\langle\lambda,\alpha'_1\rangle)M_{P_1|P}(\lambda)c_sI_P(f,\lambda))_{\mathcal A^2_P}$$
over 
\begin{itemize}
\item a set of representatives of parabolic subgroups $P$,
\item an element of the Weyl group $W(M)$,
\item an parameter $\lambda\in (\mathfrak a^G_M)^*$ 
\item a $k$-tuples of co-roots $\check{\alpha}_1,\dots,\check{\alpha}_k\in \Sigma^G_M$
\end{itemize} }

\subsection{The Orbit Method for reductive Lie groups}
Let us now describe the Orbit Method for reductive groups. The results are known
in the literature. We refer the readers to \cite{vogan} for a detailed classification of irreducible unitary representations of reductive Lie groups.

In particular, the following result are well-known.
\begin{lem}
In the matrix realization of reductive Lie groups, the dual space is identified with itself under the sesqui-linear homomorphism 
\begin{equation} X\in \mathfrak g \mapsto \langle ., X\rangle  \end{equation}
The adjoint representation and coadjoint representation become matrix conjugation 
\begin{equation} \Ad X \mapsto X(\; .\; ) X^{-1}. \end{equation}
\end{lem}

\begin{lem}
The stabilizer of a regular element is an abelian subgroup and the the corresponding orbit has a polarization that is a parabolic subgroup.
\end{lem}

\subsection{Automorphic representations via the Procedure of Quantization}
As conclusion we see that the automorphic representations could be obtained by this construction. More precisely, our main result is
\begin{thm} \label{thm51}
There is a one-to-one correspondence between the automorphic representations of $G$ and the representations obtained from the procedure of quantization based on the Geometric Langlands Correspondence.
\end{thm}
{\sc Proof.}

Let us consider an automorphic representation $T$. From the above construction, it is obtained by some induction as a sub-quotient 
$I_P(\lambda)$
of an (some kind, hololorphic) induction, see (\ref{ind})
$\mathcal A^2_P = \Ind_P^G L^2_{disc}(A_M\Gamma_M \setminus M)$
from a parabolic subgroup $P$,

If the parabolic subgroup is maximal the representation is itself a discrete series representation and it is well-known that it can be realized as an authomorphic representation. 
Let us denote by $G=K.P$ a Gauss decomposition, where $K$ is a maximal compact subgroup.

The quotient $P\setminus G$ is of complex dimension 1 and is a Riemann surface. 
After the construction, every automorphic representation via the Eisenstein series construction, is equivalent to  a sub-quotient of an induced representation from a discrete series representation of a parabolic subgroup.

Following the Geometric Langlands Correspondence \cite{kapustinwitten}, we have following results:
\begin{lem}
there is a bijection between the equivalence class of representations of the fundamental group of the Riemann surface $P\setminus G/\Gamma \cong \Gamma \cap K\setminus K_{\mathbb C}$ into the dual Langlands group ${}^LG$ and the automorphic representation of $G$ as a sub-quotient of  the induced representation (\ref{ind}). 
\end{lem}
Because of the construction of automorphic representations, there is a unique representation $\pi \in L^2_{disc}(A_M\Gamma_M \setminus M)$ such that $T = T_\pi = \ind_P^G(A_M\Gamma_M\setminus M)$ under the above bijection.

\begin{lem}
The discrete series representation $\pi\in L^2_{disc}(A_M\Gamma_M\setminus M)$ is realized by the Geometric Quantization construction in the Orbit Method.
\end{lem}

If the parabolic subgroup is not maximal, we find a tower 

of parabolic subgroup which are maximal each-into-another . 

\begin{lem}
There is a sequence of  parabolic subgroups  $P_i$,
\begin{equation} P \subset P_1 \subset \dots \subset P_k \end{equation}
such that are an maximal parabolic each-to-the next $P_i \subset P_{i+1}$ .
\end{lem}
Following the step-induction 
\begin{equation} \ind_P^GL^2(A_M\Gamma_M\backslash M) = \ind^G_{P_k} \ind^{P_k}_{P_{k-1}} \dots \ind_{P}^{P_1}L^2(A_M\Gamma_M\backslash M)  \end{equation}
In each induction step the quotient is some Riemann surface and we can use our procedure of quantization to realize the representations in the corresponding Fock spaces.
\hfill$\Box$
\begin{rem}
This is equivalent to a complete symmetry breaking in field theory and therefore we can use quantization procedure of fields, repeatedly.  
\end{rem}


\vskip 1cm
{\sc Institute of Mathematics, VAST, 18 Hoang Quoc Viet Road, Cau Giay Disctrict, 10307, Hanoi, Vietnam}\\
{\tt Email: dndiep@math.ac.vn}
\end{document}